\documentclass[a4paper,11pt]{article}

\usepackage{amssymb,amsmath,diagrams,times,theorem}
\usepackage{pstcol,pst-node,pst-grad}

\newcommand{\N}{\mathbb{N}}

\newcommand{\Oscr}{{\mathcal O}}

\newcommand{\wis}[1]{{\text{\em \usefont{OT1}{cmtt}{m}{n} #1}}}
 \newcommand{\proofhead}[1]{\par\pagebreak[1]\noindent{\bf#1.\ }}
 \newcommand{\pf}{\proofhead{Proof}}
  \newcommand{\qed}{{\unskip\nolinebreak[1]\hspace{1.5em}\mbox{}\nolinebreak
    \hfill$\Box$\parfillskip=0pt\finalhyphendemerits=0\par\pagebreak[1]}}

\newenvironment{proof}{\pf}{\qed}

\newtheorem{theorem}{Theorem}
\newtheorem{proposition}{Proposition}

{\theorembodyfont{\rmfamily} 
\newtheorem{example}{Example}

\newtheorem{definition}{Definition} }

\title{A non-commutative topology on $\wis{rep}$~A}

\author{Lieven Le Bruyn \\
Departement Wiskunde en Informatica \\ Universiteit Antwerpen  \\
B-2020 Antwerp (Belgium) \\
lieven.lebruyn@ua.ac.be}

\date{}

\begin{document}

\maketitle

\begin{abstract}
We extend the Zariski topology on $\wis{simp}~A$, the set of all simple finite dimensional representations of $A$, to a non-commutative topology (in the sense of Fred Van Oystaeyen) on $\wis{rep}~A$, the set of all finite dimensional representations of $A$, using Jordan-H\"older filtrations. The non-commutativity of the topology is enforced by the order of the composition factors.
\end{abstract}

All algebras will be affine associative $\Bbbk$-algebras with unit over an algebraically closed field $\Bbbk$. The {\em non-commutative affine 'scheme'} associated to an algebra $A$ is, as a set,  the disjoint union
\[
\wis{rep}~A = \bigsqcup_n~\wis{rep}_n~A \]
where $\wis{rep}_n~A$ is the (commutative) affine scheme of $n$-dimensional representations of $A$.  In this note we will equip $\wis{rep}~A$ with a non-commutative topology in the sense of Fred Van Oystaeyen \cite[\S 7.2]{FVOAGAA} (or, more precisely, a slight generalization of it).

Here is the main idea. The twosided prime ideal spectrum $\wis{spec}~A$ is an (ordinary) topological space via the Zariski topology, see for example \cite{FVO444} or \cite[\S II.6]{ProcesiBook}. Hence, the subset $\wis{simp}~A$ of all simple finite dimensional $A$-representations can be equipped with the induced topology. This topology can then be extended to a non-commutative topology on $\wis{rep}~A$ using Jordan-H\"older filtrations. The non-commutative nature of the topology is enforced by the order of the composition factors. 

We give a few examples, connect this notion with that of Reineke's composition monoid and remark on the difference between quotient varieties and moduli spaces from the perspective of non-commutative topology. Finally, we note that this construction can be generalized verbatim to any Artinian Abelian category as soon as we have a topology on the set of simple objects.

\section{The Zariski topology on $\wis{simp}~A$.}

Recall that a prime ideal $P$ of $A$ is a twosided ideal satisfying the property that if $I.J \subset P$ then $I \subset P$ or $J \subset P$ for any pair of twosided ideals $I,J$ of $A$. The {\em prime spectrum} $\wis{spec}~A$ is the set of all twosided prime ideals of $A$. The {\em Zariski topology} on $\wis{spec}~A$ has as its closed subsets
\[
\mathbb{V}(S) = \{ P \in \wis{spec}~A~|~S \subset P \} \]
where $S$ varies over all subsets of $A$, see for example \cite[Prop. II.6.2]{ProcesiBook}. Note that an algebra morphism $\phi~:~A \rTo B$ does {\em not} necessarily induce a continuous map $\phi^*~:~\wis{spec}~B \rTo \wis{spec}~A$ but is does so in the case $\phi$ is a {\em central extension} in the sense of \cite[\S II.6]{ProcesiBook}.

If $M \in \wis{rep}_n~A$ is a simple $n$-dimensional representation, there is a defining epimorphism $\psi_M~:~A \rOnto M_n(\Bbbk)$ and the kernel of this morphism $\wis{ker}~\psi_M$ is a twosided maximal (hence prime) ideal of $A$. We define the Zariski topology on the set of all simple finite dimensional representations $\wis{simp}~A$ by taking as its closed subsets
\[
\mathbb{V}(S) = \{ M \in \wis{simp}~A~|~S \subset \wis{ker}~\psi_M \} \]
Again, one should be careful that whereas an algebra map $\phi~:~A \rTo B$ induces a map $\phi^*~:~\wis{rep}~B \rTo \wis{rep}~A$ it does {\em not} in general map $\wis{simp}~B$ to $\wis{simp}~A$ (unless $\phi$ is a central extension).

With $\mathcal{L}_A$ we will denote the set of all open subsets of $\wis{simp}~A$. $\mathcal{L}_A$ will be the set of {\em letters} on which to base our non-commutative topology.

\section{Non-commutative topologies (and generalizations).}

In \cite[Chp. 7]{FVOAGAA} Fred Van Oystaeyen defined {\em non-commutative topologies} which are generalizations of usual topologies in which it is no longer  true that $A \cap A$ is equal to $A$ for an open set $A$. In order to keep  dichotomies of possible definitions to a minimum he imposed left-right symmetric conditions on the definition. However, for applications to representation theory it seems that the most natural non-commutative topologies are truly one-sided. For this reason we take some time to generalize some definitions and results of \cite[Chp. 7]{FVOAGAA}.

We fix a partially ordered set $(\Lambda,\leq)$ with a unique minimal element $0$ and a unique maximal element $1$, equipped with two operations $\wedge$ and $\vee$. With $i_{\Lambda}$ we will denote the set of all {\em idempotent elements} of $\Lambda$, that is, those $x \in \Lambda$ such that $x \wedge x = x$. A {\em finite global cover} is a finite subset $\{ \lambda_1,\hdots,\lambda_n \}$ such that $1 = \lambda_1 \vee \hdots \vee \lambda_n$. In the table below we have listed the conditions for a (one-sided) non-commutative topology. Note that some requirements are less essential than others. For example, the covering condition (A10) is only needed if we want to fit non-commutative topologies in the framework of non-commutative Grothendieck topologies \cite{FVOAGAA} and the weak modularity condition (A9) is not required if every basic open is $\vee$-idempotent (as is the case in most examples). 

\begin{Rotateright}
\parbox{18cm}{
\[
\begin{array}{l|ccc}
(A1) & x \wedge y \leq x & & x \wedge y \leq y \\ \\
\hline \\
(A2) & x \wedge 1 = x & & 1 \wedge x = x \\
& x \wedge 0 = 0 & & 0 \wedge x = 0 \\ \\
\hline \\
(A3) & & (x \wedge y) \wedge z = x \wedge (y \wedge z) = x \wedge y \wedge z & \\ \\
\hline \\
(A4) & x \leq y \Rightarrow z \wedge x \leq z \wedge y & & x \leq y \Rightarrow x \wedge z \leq y \wedge z \\ \\
\hline \\
(A5) & x \leq x \vee y & & y \leq x \vee y \\ \\
\hline \\
(A6) & x \vee 1 = 1 & & 1 \vee x = 1 \\
& x \vee 0 = x & & 0 \vee x = x \\ \\
\hline \\
(A7) & & (x \vee y) \vee z = x \vee (y \vee z) = x \vee y \vee z \\ \\
\hline \\
(A8) & x \leq y \Rightarrow x \vee z \leq y \vee z & & x \leq y \Rightarrow z \vee x \leq z \vee y \\ \\
\hline \\
(A9) & a \vee (a \wedge b) \leq (a \vee a) \wedge b & & a \vee (b \wedge a) \leq (a \vee b) \wedge a \\ \\
\hline \\
(A10) & x = (x \wedge \lambda_1) \vee \hdots \vee (x \wedge \lambda_n) & & 
x = (\lambda_1 \wedge x) \vee \hdots \vee (\lambda_n \wedge x)
\end{array}
\] }
\end{Rotateright} 

\begin{definition} Let $(\Lambda,\leq)$ be a partially ordered set with minimal and maximal element $0$ and $1$ and operations $\wedge$ and $\vee$. Then,

$\Lambda$ is said to be a {\em left non-commutative topology} if and only if the left and middle column conditions of (A1)-(A10) are valid for all $x,y,z \in \Lambda$, all $a,b \in i_{\Lambda}$ with $a \leq b$ and all finite global covers $\{ \lambda_1,\hdots,\lambda_n \}$.

$\Lambda$ is said to be a {\em right non-commutative topology} if and only if the middle and right column conditions of (A1)-(A10) are valid for all $x,y,z \in \Lambda$, all $a,b \in i_{\Lambda}$ with $a \leq b$ and all finite global covers $\{ \lambda_1,\hdots,\lambda_n \}$.

$\Lambda$ is said to be a {\em non-commutative topology} if and only if the conditions (A1)-(A10) are valid for all $x,y,z \in \Lambda$, all $a,b \in i_{\Lambda}$ with $a \leq b$ and all finite global covers $\{ \lambda_1,\hdots,\lambda_n \}$.
\end{definition}

There are at least two ways of building a genuine non-commutative topology out of these sets of basic opens. We briefly sketch the procedures here and refer to the forthcoming monograph \cite{FredFunctor} for details in the symmetric case (the one-sided versions present no real problems).

Let $T(\Lambda)$ be the set of all finite $(\wedge,\vee)$-words in the {\em contractible} idempotent elements $i_{\Lambda}$ (that is, $\lambda \in i_{\Lambda}$ such that for all $\lambda_1,\lambda_2$ with $\lambda \leq \lambda_1 \vee \lambda_2$ we have that $\lambda = (\lambda \wedge \lambda_1) \vee (\lambda \wedge \lambda_2)$). If $\Lambda$ is a (left,right) non-commutative topology, then so is $T(\Lambda)$. The {\em $\vee$-complete topology of virtual opens} $T'(\Lambda)$ is then the set of all $(\wedge,\vee)$-words in the contractible idempotents of finite length in $\wedge$ (but not necessarily of finite length in $\vee$). This non-commutative topology has properties very similar to that of an ordinary topology and, in fact, has associated to it a {\em commutative shadow}.

The second construction, leading to the {\em pattern topology}, starts with the equivalence classes of {\em directed systems} $S \subset \Lambda$ (that is, if for all $x,y \in S$ there is a $z \in S$ such that $z \leq x$ and $z \leq y$) and where the equivalence relation $S \sim S'$ is defined by
\[
\begin{cases}
\forall a \in S, \exists a' \in S, a' \leq a~\text{and}~b \leq a' \leq b'~\text{for some}~b,b' \in S' \\
\forall b \in S', \exists b' \in S', b' \leq b~\text{and}~a \leq b' \leq a'~\text{for some}~a,a' \in S
\end{cases}
\]
One can extend the $\wedge,\vee$ operations on $\Lambda$ to the equivalence classes $C(\Lambda) = \{ [S]~|~S $ directed $\}$ in the obvious way such that also $C(\Lambda)$ is a (left,right) non-commutative topology. A directed set $S \subset \Lambda$ is said to be {\em idempotent} if for all $a \in S$, there is an $a' \in S \cap i_{\Lambda}$ such that $a' \leq a$. If $S$ is idempotent then $[S] \in i_{C(\Lambda)}$ and those idempotents will be called {\em strong idempotents}. The pattern topology $\Pi(\Lambda)$ is the (left,right) non-commutative topology of finite $(\wedge,\vee)$-words in the strong idempotents of $C(\Lambda)$. A directed system $[S]$ is called a {\em point} iff $[S] \leq \vee [S_{\alpha}]$ implies that $[S] \leq [S_{\alpha}]$ for some $\alpha$.

\section{The basic opens.}

For an $n$-dimensional representation $M$ of $A$ we call a finite filtration of length $u$
\[
\mathcal{F}^u~:~0 = M_0 \subset M_1 \subset \hdots \subset M_u = M \]
of $A$-representations a {\em Jordan-H\"older filtration} if the successive quotients
\[
\mathcal{F}_i = \frac{M_i}{M_{i-1}} \]
are simple $A$-representations. Recall that $\mathcal{L}_A$ is the set of all open subsets $V$ of $\wis{simp}~A$. With $\mathbb{W}_A$ we denote the non-commutative words in these letters
\[
\mathbb{W}_A = \{ V_1 \hdots V_k~|~V_i \in \mathcal{L}_A, k \in \N \} \]
For a given word $w = V_1 V_2 \hdots V_k \in \mathbb{W}_A$ we define the {\em left basic open set}
\[
\Oscr_w^l = \{ M \in \wis{rep}~A~|~\exists \mathcal{F}^u~\text{Jordan-H\"older filtration on $M$ such that}~\mathcal{F}_i \in V_i \}
\]
and the {\em right basic open set}
\[
\Oscr_w^r = \{ M \in \wis{rep}~A~|~\exists \mathcal{F}^u~\text{Jordan-H\"older filtration on $M$ such that}~\mathcal{F}_{u-i} \in V_{k-i} \}
\]
Finally, to make these definitions symmetric we define the {\em basic open set}
\[
\Oscr_w = \{ M \in \wis{rep}~A~|~\exists \mathcal{F}^u~\text{Jordan-H\"older filtration on $M$ such that}~\mathcal{F}_{i_j} \in V_{j} \]
\[
\qquad~\quad~\text{for some}~1 \leq i_1 < i_2 < \hdots < i_k \leq u~ \}
\]
Clearly, $\Oscr_w^l$ consists of those representations having prescribed  bottom structure, whereas $\Oscr_w^r$ consists of those with prescribed top structure. In order to avoid three sets of definitions we will denote from now on $\Oscr_w^{\bullet}$ whenever we mean $\bullet \in \{ l,r,\emptyset \}$.

If $w=L_1\hdots L_k$ and $w'= M_1 \hdots M_l$, we will denote with $w \cup w'$ the {\em multi-set}
$\{ N_1,\hdots,N_m \}$ where each $N_i$ is one of $L_j,M_j$ and $N_i$ occurs in $w \cup w'$ as many times as its maximum number of factors in $w$ or $w'$. With $\wis{rep}(w \cup w')$ we denote the subset of $\wis{rep}~A$ consisting of the representations of $M$ having a Jordan-H\"older filtration having factor-multi-set containing $w \cup w'$. For any triple of words $w,w'$ and $w"$ we denote $\Oscr_{w"}^{\bullet}(w \cup w') = \Oscr^{\bullet}_{w"} \cap \wis{rep}(w \cup w')$.

We define an equivalence relation on the basic open sets by
\[
\Oscr^{\bullet}_w \approx \Oscr^{\bullet}_{w'} \qquad \Leftrightarrow \qquad
\Oscr^{\bullet}_w(w \cup w') = \Oscr^{\bullet}_{w'}(w \cup w') \]
The reason for this definition is that the condition of $M \in \Oscr^{\bullet}_w$ is void if $M$ does not have enough Jordan-H\"older components to get all factors of $w$ which makes it impossible to define equality of basic open sets defined by different words.

We can now define the partially ordered sets $\Lambda^{\bullet}_A$ as consisting of all basic open subsets $\Oscr_w^{\bullet}$ of $\wis{rep}~A$. The partial ordering $\leq$ is induced by set-theoretic inclusion modulo equivalence, that is,
\[
\Oscr^{\bullet}_w \leq \Oscr^{\bullet}_{w'} \qquad \Leftrightarrow \qquad \Oscr^{\bullet}_w(w \cup w') \subseteq \Oscr^{\bullet}_{w'}(w \cup w') \]
As a consequence, equality $=$ in the set $\Lambda_A^{\bullet}$ coincides with equivalence $\approx$. Observe that these partially ordered sets have a unique minimal and a unique maximal element (upto equivalence)
\[
0 = \emptyset = \Oscr_{\emptyset}^{\bullet} \qquad \text{and} \qquad 1 = \wis{rep}~A = \Oscr_{\wis{simp}~A}^{\bullet}
\]
The operations $\vee$ and $\wedge$ are defined as follows : $\vee$ is induced by ordinary set-theoretic union and $\wedge$ is induced by concatenation of words, that is
\[
\Oscr_w^{\bullet} \wedge \Oscr_{w'}^{\bullet} \approx \Oscr_{ww'}^{\bullet}  \]

\begin{theorem} With notations as before :
\begin{itemize}
\item{$(\Lambda_A^l,\leq,\approx,0,1,\vee,\wedge)$ is a left non-commutative topology on $\wis{rep}~A$.}
\item{$(\Lambda_A^r,\leq,\approx,0,1,\vee,\wedge)$ is a right non-commutative topology on $\wis{rep}~A$.}
\end{itemize}
\end{theorem}

\begin{proof} The tedious verification is left to the reader. Here, we only stress the importance of the equivalence relation for example in verifying $x \wedge 1 = x$. So, let $w=L_1\hdots L_k$ then
\[
\Oscr_w^{l} \wedge 1 = \Oscr_{L_1\hdots L_k \wis{simp} A}^{l} \subset \Oscr_{w}^{l} \]
and this inclusion is proper (look at elements in $\Oscr_w^{l}$ having exactly $k$ composition factors). However, as soon as the representation has $k+1$ composition factors, it is contained in the left hand side whence $\Oscr_w^{l} \wedge 1 \approx \Oscr_w^{l}$. A similar argument is needed in the covering condition.
\end{proof}

Note however that $(\Lambda_A,\leq,\approx,0,1,\vee,\wedge)$ is not necessarily a non-commutative topology : the problematic conditions are $\Oscr_w \wedge 1 = \Oscr_w = 1 \wedge \Oscr_w$ and the covering condition. The reason is that for $w=L_1\hdots L_k$ as before and $M \in \Oscr_w$  having $> k$ factors, it may happen that the last factor is the one in $L_k$ leaving no room for a successive factor in $\wis{simp}~A$ (whence $\Oscr_w \cap 1$ is not equivalent to $\Oscr_w$).

\begin{example} Let $A$ be a finite dimensional algebra, then $A$ has a finite number of simple representations $\wis{simp}~A = \{ S_1,\hdots,S_n \}$ and the Zariski topology is the discrete topology. If for some $1 \leq i ,j \leq n$ we have that
\[
Ext^1_A(S_i,S_j) = 0 \qquad \text{and} \qquad Ext^1_A(S_j,S_i) \not= 0 \]
then $\Lambda_A^l$ is a genuinely non-commutative topology, for example
\[
\Oscr_{S_i}^l \wedge \Oscr_{S_j}^l = \Oscr_{S_iS_j}^l \not=  \Oscr_{S_jS_i}^l = \Oscr_{S_j}^l \wedge  \Oscr_{S_i}^l \]
as a non-trivial extension $0 \rTo S_i \rTo X \rTo S_j \rTo 0$ belongs to $\Oscr_{S_iS_j}^l(S_iS_j \cup S_jS_i)$ but not to $\Oscr_{S_jS_i}^l(S_iS_j \cup S_jS_i)$.
\end{example}

\section{Reineke's mon(str)oid.}

When $A$ is the path algebra of a quiver without oriented cycles we can generalize the foregoing example and connect the previous definitions to the {\em composition monoid} introduced and studied by Markus Reineke in \cite{Reineke}.

Let $Q$ be a quiver without oriented cycles, then its path algebra $A = \Bbbk Q$ is finite dimensional hereditary with all simple representations one-dimensional and in one-to-one correspondence with the vertices of $Q$. For every dimension $n$ we have that
\[
\wis{rep}_n~A = \bigsqcup_{|\alpha | = n} GL_n \times^{GL(\alpha)} \wis{rep}_{\alpha}~Q \]
where $\alpha$ runs over all dimension vectors of total dimension $n$ and where $\wis{rep}_{\alpha}~Q$ is the affine space of all $\alpha$-dimensional representations of the quiver $Q$ with base-change group action by $GL(\alpha)$.

The {\em Reineke monstroid} $\mathcal{M}(Q)$ has as its elements the set of all irreducible closed $GL(\alpha)$-stable subvarieties of $\wis{rep}_{\alpha}~Q$ for all dimension vectors $\alpha$, equipped with a product
\[
\mathcal{A} \ast \mathcal{B} = \{ X \in \wis{rep}_{\alpha+\beta}~Q~|~\text{there is an exact sequence} \]
\[
0 \rTo M \rTo X \rTo N \rTo 0~\quad M \in \mathcal{A},N \in \mathcal{B}~\} \]
if $\mathcal{A}$ (resp. $\mathcal{B}$) is an element of $\mathcal{M}(Q)$ contained in $\wis{rep}_{\alpha}~Q$ (resp. in $\wis{rep}_{\beta}~Q$). It is proved in \cite[lemma 2.2]{Reineke} that $\mathcal{A} \ast \mathcal{B}$ is again an element of $\mathcal{M}(Q)$. This defines a monoid structure on $\mathcal{M}(Q)$ which is too unwieldy to study directly. Observe that we changed the order of the terms wrt. the definition given in \cite{Reineke}. That is, we will work with the {\em opposite} monoid of \cite{Reineke}.

On the other hand, the {\em Reineke composition monoid} is very tractable. It is the submonoid $\mathcal{C}(Q)$ of $\mathcal{M}(Q)$ generated by the vertex-representation spaces $R_i = \wis{rep}_{\delta_i}~Q$. These generators satisfy specific commutation relations which can be read off from the quiver structure, see \cite[\S 5]{Reineke}. For example, if there are no arrows between $v_i$ and $v_j$ then
\[
R_i \ast R_j = R_j \ast R_i \]
and if there are no arrows from $v_i$ to $v_j$ but $n$ arrows from $v_j$ to $v_i$, then
\[
\begin{cases}
R_i^{\ast (n+1)} \ast R_j = R_i^{\ast n} \ast R_j \ast R_i \\
R_i \ast R_j^{\ast (n+1)} = R_j \ast R_i \ast R_j^{\ast n}
\end{cases}
\]
For more details on the structure of $\mathcal{C}(Q)$ we refer to \cite[\S 5]{Reineke}.

There is a relation between $\mathcal{C}(Q)$ and the left- and right- non-commutative topologies $\Lambda_A^l$ and $\Lambda_A^r$. Because the Zariski topology on $\wis{simp}~A$ is the discrete topology on the set $\{ S_1,\hdots,S_k \}$ of vertex simples, it is important to understand $\Oscr_w^r$ where $w$ is a word in the $S_i$, say
$w = S_{i_1}S_{i_2}\hdots S_{i_u}$. In fact, we could have based our definition of a one-sided non-commutative topology on the set $\mathcal{L}_A$ of {\em irreducible} open subsets of $\wis{simp}~A$ and then these basic opens would be all. If $\mathcal{C}$ is a $GL(\alpha)$-stable subset of $\wis{rep}_{\alpha}~Q$ with $| \alpha | = n$, we will denote the subset $GL_n \times^{GL(\alpha)} \mathcal{C}$ of $\wis{rep}_n~A$ by $\tilde{\mathcal{C}}$.

\begin{proposition}
\[
\Oscr^l_w = \bigcup_{w'} \tilde{\mathcal{A}}_{w'} \qquad \text{resp.} \qquad \Oscr^r_w = \bigcup_{w'} \tilde{\mathcal{A}}_{w'} \]
where $\mathcal{A}_{w'}$ is a $\ast$-word in the generators $R_i$ of the composition monoid such that $w'$ can be rewritten (using the relations in $\mathcal{C}(Q)$) in the form
\[
w' = R_{i_1} \ast R_{i_2} \ast \hdots \ast R_{i_u} \ast w" \quad \text{resp.} \quad 
w' = w" \ast R_{i_1} \ast R_{i_2} \ast \hdots \ast R_{i_u}  \]
for another $\ast$-word $w"$. 
\end{proposition}

Also, the equivalence relation introduced before can be expressed in terms of $\mathcal{C}(Q)$. If $w = S_{i_1}S_{i_2} \hdots S_{i_u}$ and $w' = S_{j_1}S_{j_2}\hdots S_{j_v}$ such that $w \cup w' = \{ S_{k_1},\hdots,S_{k_w} \}$, then

\begin{proposition} $\Oscr^l_w \approx \Oscr^l_{w'}$ if and only if every $\ast$-word $v = R_{a_1} \ast \hdots \ast R_{a_z}$ containing in it distinct factors $R_{k_1},\hdots,R_{k_w}$ which can be brought in $\mathcal{C}(Q)$ in the form
\[
v = R_{i_1} \ast \hdots \ast R_{i_u} \ast v' \]
can also be written in the form
\[
v = R_{j_1} \ast \hdots \ast R_{j_v} \ast v" \]
(and conversely). A similar result describes $\Oscr^r_w \approx \Oscr^r_{w'}$.
\end{proposition}

In particular, in this setting there will be hardly any {\em idempotent} basic opens (that is, satisfying
$\Oscr^r_w \wedge \Oscr^r_w \approx \Oscr^r_w$). Clearly, if $\{ S_{e_1},\hdots,S_{e_a} \}$ are simples such that the quiver restricted to $\{ v_{e_1},\hdots,v_{e_a} \}$ has no arrows, then any word $w$ in the $S_{e_j}$ gives an idempotent $\Oscr^r_w$. In the following section we will give an example where {\em every} basic open is idempotent and hence we get a commutative topology.

\section{The commutative case.}

If $A$ is a commutative affine $\Bbbk$-algebra, then any simple representation is one-dimensional, $\wis{simp}~A = X_A$ the affine (commutative) variety corresponding to $A$ and the Zariski topologies on both sets coincide. Still, one can define the non-commutative topologies on $\wis{rep}~A$. However,

\begin{proposition} If $A$ is a commutative affine $\Bbbk$-algebra, then  both $\Lambda_A^l$ and $\Lambda^r_A$ are commutative topologies. That is, for all words $w$ and $w'$ in $\mathcal{L}_A$ we have
\[
\Oscr_w^l \wedge \Oscr^l_{w'} \approx \Oscr^l_{w'} \wedge \Oscr^l_{w} \qquad \text{and} \qquad 
\Oscr_w^r \wedge \Oscr^r_{w'} \approx \Oscr^r_{w'} \wedge \Oscr^r_{w}
\]
\end{proposition}

\begin{proof}
We claim that every basic open $\Oscr^l_w$ is idempotent. Observe that all simple $A$-representations are one-dimensional and that there are only self-extensions of those, that is, if $S$ and $T$ are non-isomorphic simples, then $Ext^1_A(S,T)=0=Ext^1_A(T,S)$. However, there are self-extensions with the dimension of $Ext^1_A(S,S)$ being equal to the dimension of the tangent space at $X_A$ in the point corresponding to $S$. As a consequence we have for any Zariski open subsets $U$ and $V$ of $X_A$ that
\[
\Oscr^l_{UV} = \Oscr^l_{VU} \]
as we can change the order of the filtration factors (a representation $M$ is the direct sum of submodules $M_1 \oplus \hdots \oplus M_s$ with each $M_i$ concentrated in a single simple $S_i$ and we can add the successive $S_i$ factors of $M$ at any wanted place in the filtration sequence). Hence, for every word $w$ we have that
\[
\Oscr^l_w \approx \Oscr^l_w \wedge \Oscr_w^l \]
and also for any pair of words $w$ and $w'$ we have that
\[
\Oscr_w^l \wedge \Oscr_{w'}^l = \Oscr^l_{ww'} = \Oscr_{w'w}^l = \Oscr_{w'}^l \wedge \Oscr^l_{w}
\]
Observe that in \cite{FVOAGAA} it is proved that a non-commutative topology in which every basic open is idempotent is commutative. We cannot use this here as the proof of that result  uses both the left- and right- conditions. However, we are dealing here with a very simple example.
\end{proof}

\section{Quotient varieties versus moduli spaces.}

Having defined a one-sided non-commutative topology on $\wis{rep}~A$ we can ask about the induced topology on the quotient variety $\wis{iss}~A$ of all isomorphism classes of semi-simple $A$-representations or on the moduli space $\wis{moduli}_{\theta}~A$ with respect to a certain stability structure $\theta$, cfr. \cite{Rudakov}. Experience tells us that it is a lot easier to work with quotient varieties than with moduli spaces and non-commutative topology may give a partial explanation for this.

Indeed, as the points of $\wis{iss}~A$ are semi-simple representations, it is clear that the induced non-commutative topology on $\wis{iss}~A$ is in fact commutative. However, as the points of $\wis{moduli}_{\theta}~A$ correspond to isomorphism classes of direct sums of stable representations (not simples!), the induced non-commutative topology on $\wis{moduli}_{\theta}~A$ will in general remain non-commutative. Still, in nice examples, such as representations of quivers, one can 
define another non-commutative topology on $\wis{moduli}_{\theta}~A$ which does become commutative. Use universal localization to cover $\wis{moduli}_{\theta}~A$ by opens isomorphic to $\wis{iss}~A_{\Sigma}$ for some families $\Sigma$ of maps between projectives and equip $\wis{moduli}_{\theta}~A$ with a non-commutative topology (which then will be commutative!) obtained by gluing the induced non-commutative topologies on the $\wis{rep}~A_{\Sigma}$.

\section{Generalizations.}

It should be evident that our construction can be carried out verbatim in the setting of any Artinian Abelian category (that is, an Abelian category having Jordan-H\"older sequences) as soon as we have a natural topology on the set of simple objects. In fact, the same procedure can be applied when we have a left (or right) non-commutative topology on the simples.

In fact, the construction may even be useful in Abelian categories in which every object is filtered by special objects on which we can define a (one-sided) (non-commutative) topology.

\end{document}